# Cauchy-like Criterion for Differentiability of Functions of Several Variables

*Yurii V. Mukhin[†] and Nataliya D. Kundikova*

South Ural State University: *Chelyabinsk, Russian Federation*
Institute of Electrophysics, UB RAS: *Ekaterinburg, Russian Federation*

## *Abstract*

In this paper, several differentiability criteria for real functions of multiple variables in n-dimensional Euclidean space are considered. Simple and easy-to-use Cauchy-like criterion is formulated and proven. Relaxed sufficient conditions for differentiability that do not require continuity of all partial derivatives are suggested. Generalization of the Cauchy-like criterion for functions on cross products of normed vector spaces (not necessarily Banach spaces) is discussed. The results of this study can be used in systems analysis, linear programming, optimization methods, functional analysis, topology and convex analysis.

## *Used Notations*

Without loss of generality, we consider differentiability of a real function f($\vec{x}$) at a point $\vec{x}=0$, where $\vec{x}$ is a vector in n-dimensional Euclidian space: $\vec{x} \in \mathbb{R}^n$ (f: U → $\mathbb{R}$, where U is some vicinity of the origin, U ⊂ $\mathbb{R}^n$). Also, without loss of generality, real function f($\vec{x}$) is assumed to be equal to zero at the origin $\vec{x}=0$, f(0)=0. Different vectors $\vec{x}_j$ are distinguished by different subscript indexes, while their components (space coordinates) are marked with different superscript indexes: $\vec{x}_j = \{x_j^1,... x_j^i,... x_j^n\}$.

In our consideration, we use asymptotic notations of big $\mathcal{O}(\varrho)$ and little $o(\varrho)$ that are noted as follows: $\mathcal{O}(\varrho)/\varrho$ is a bounded function when $\varrho \to 0$, while $[o(\varrho)/\varrho] \to 0$ when $\varrho \to 0$. Next, we introduce and denote partial values of any function f($\vec{x}$) as following: $f^i = f^i(\vec{x}) = f(x^i \cdot \vec{e}_i) = f(\{0,... x^i,... 0\})$, where $\vec{e}_i$ is the unit vector for the i-th coordinate axis of the vector space $\mathbb{R}^n$.

## *Introduction*

Any criterion is just the same original definition or statement only rewritten in different terms and notations. Different criteria may not only differ by terms but also by difficulty of their use. In fact, the general opinion is that there is no good or easy-to-use differentiability criterion for multivariable functions. In our paper we show otherwise.

Differentiability of multivariable functions is discussed in most of the textbooks on Calculus [1, 2]. In this section, we briefly repeat definitions and general statements with some common examples.

Differentiability for a function of n variables is defined in the following way. A function f($\vec{x}$) is said to be differentiable at the point $\vec{x}=0$ if there is a finite vector $\vec{A} \in \mathbb{R}^n$ such that:

$$f(\vec{x}) = \vec{A} \cdot \vec{x} + o(\varrho), \text{ where } \varrho = |\vec{x}| = \sqrt{\sum_{i=1}^{n}(x^i)^2} \qquad (1)$$

The definition (1) does not provide or require any information about the quality of the function in the vicinity of the origin $\vec{x}=0$. In fact, f($\vec{x}$) can be discontinuous in any close proximity of the origin. For example, let us consider a function of two real variables G(*x*,*y*) = *xy*·D(*xy*) where D(*t*) is the *Dirichlet* function of *t* ∈ $\mathbb{R}$. G(*x*,*y*) is differentiable at the origin and is discontinuous elsewhere but *x* or *y* axes.

Still, some conclusions can be made from the definition (1). Since the scalar product of two vectors in (1) is an invariant with respect to rotations of the Cartesian coordinate system, an immediate consequence of (1) is that any partial derivative exists in any direction at the origin and it is equal to $\vec{A} \cdot \vec{e}$, where $\vec{e}$ – is a unit vector in that direction. Also, the said derivative is a continuous function of the polar or spherical angles. However, these necessary conditions are not sufficient conditions for differentiability. Let us consider an example of a function of two variables:

$$g(x,y) = \frac{x^2 y}{x^2 + y^2} \text{ when } (x^2+y^2) \neq 0; \text{ and } g(0,0) = 0. \qquad (2)$$

The plot of the function g(*x*,*y*) is shown in the figure 1. This function is not differentiable at the origin. In the polar coordinates, the function is given as: g($\varrho\cos\varphi$, $\varrho\sin\varphi$) = $\varrho\cos^2\varphi \sin\varphi$. It is a linear function of $\varrho$. Thus, for any direction of given azimuthal angle $\varphi$, the derivative along the direction exists at the origin, is equal to $\cos^2\varphi \sin\varphi$, and is a continuous function of $\varphi$. Interestingly, the bijective function (2) transforms any straight line that goes through the origin in the (*x*,*y*)-plane onto a straight line that belongs entirely in the surface of $z = g(x,y)$ in three dimensional space. These straight lines are tangent lines to the surface $z = g(x,y)$ at the origin, but they do not form a common plane. That is the reason why this surface is not smooth at the origin. Thus, we have come to formulation of the geometrical criterion for differentiability.

---

[†] mukhinyv@gmail.com; mukhinyv@susu.ru; Google Scholar

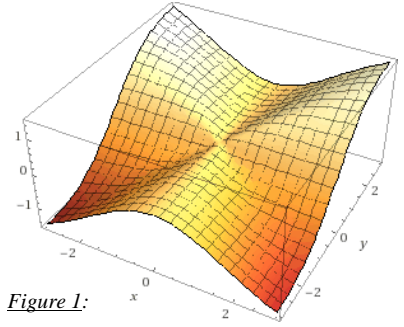

*Figure 1*:

## *Geometrical criterion for differentiability*

Geometrical criterion is associated with the construction of n-dimensional hyperplane that is a tangent plane to the surface $z = f(\vec{x})$ in the (n+1)-dimensional space $\{z, \mathbb{R}^n\}$. The tangent plane to the surface is such a plane that when moving on the surface towards a touching point, the distance to the plane vanishes faster (like a little $o$) than the distance from the touching point to the projection of the moving point onto the plane.

*Geo-Criterion*: If in (n+1)-dimensional Euclidean space $\{z, \mathbb{R}^n\}$, there is a hyperplane that is not orthogonal to $\mathbb{R}^n$ and is tangent at the origin to the surface $z = f(\vec{x})$ (where $\vec{x} \in \mathbb{R}^n$) then and only then $f(\vec{x})$ is differentiable at the origin.

The non-orthogonality requirement is necessary for the vector $\vec{A}$ in the equation (1) to have finite components. In other words, this necessary condition is implicitly included in the original definition (1), while it has to be explicitly stated in the formulation of the Geo-criterion.

The proof of the Geo-criterion is somewhat trivial (see for example [2]) and is based on the following. When a point $\{f(\vec{x}), \vec{x}\}$ is moving in the surface $z = f(\vec{x})$ towards the touching point (the origin), the distance from the touching point to the projection of the moving point $\{f(\vec{x}), \vec{x}\}$ onto the hyperplane and the $\varrho = |\vec{x}|$ are mutually big $\mathcal{O}$ asymptotic functions. We show the proof in the Appendix A at the end of our paper.

It is an important question whether a criterion is easy or difficult to use in practice. The Geo-criterion is not very easy to use as well as the original differentiability definition (1). In that respect, it is interesting to consider another criterion – the Cauchy criterion.

## *Cauchy criterion for differentiability*

The Cauchy criterion for differentiability was described by A.M. Macbeath in 1956 [3]. Let us consider (n+1) small, different and independent vectors $\vec{x}_j$ ($\vec{x}_j \in \mathbb{R}^n$) and let a function $f(\vec{x})$ be differentiable at the origin. As before, we assume that $f(0)=0$. Since the condition (1) is true for all of the vectors $\vec{x}_j$, we can write that in a form of a system of (n+1) equations as shown in (3).

$$\begin{pmatrix} f(\vec{X}_1) & x_1^1 & ... & x_1^i & ... & x_1^n \\ .... & ... & ... & ... & ... & ... \\ f(\vec{X}_j) & x_j^1 & ... & x_j^i & ... & x_j^n \\ .... & ... & ... & ... & ... & ... \\ f(\vec{X}_{n+1}) & x_{n+1}^1 & ... & x_{n+1}^i & ... & x_{n+1}^n \end{pmatrix} \begin{pmatrix} -1 \\ \vec{A} \end{pmatrix} = \begin{pmatrix} o(\varrho_1) \\ ... \\ o(\varrho_j) \\ ... \\ o(\varrho_n) \end{pmatrix}, \quad \text{where } \varrho_j = \left|\vec{X}_j\right| \qquad (3)$$

Let us denote the determinant of the (n+1)×(n+1) matrix in the left side of (3) as $\mathcal{D}(f: \vec{x}, n)$. If $f(\vec{x})$ is a linear function then the right side of (3) is identically equal to zero. Therefore, the system of the (n+1) algebraic equations with respect to $(-1, \vec{A})$ becomes uniform. That means that $\mathcal{D}(f: \vec{x}, n)$ must be equal to zero for any arbitrary set of $\vec{x}_j$. In general, when $f(\vec{x})$ is not a linear function, and when all of the $\vec{x}_j$ are close to the origin, the $\mathcal{D}(f: \vec{x}, n)$ must behave like a product of all the $o(\varrho_j)$ or like little $o$ of the product of all the $\varrho_j$: $o(\prod_j |\vec{x}_j|)$. It turns out that this necessary condition is also a sufficient condition for differentiability.

Finally, the formulation of the Cauchy criterion is as follows.

*Cauchy criterion*: For an arbitrary set of (n+1) vectors $\vec{x}_j \in \mathbb{R}^n$ and a function $f(\vec{x})$ that is equal to zero at the origin ($f(0)=0$), if $\mathcal{D}(f: \vec{x}, n)$ satisfies the condition (4) then and only then $f(\vec{x})$ is differentiable at the origin.

$$|\mathcal{D}(f: \vec{x}, n)| = o(\prod_{j=1}^{n+1} |\vec{x}_j|), \quad \text{when } \max_j\{|\vec{x}_j|\} \to 0 \qquad (4)$$

That is essential that the product of all $\varrho_j$ is there at the right side of the equation (4) and not something like $(\max_j\{\varrho_j\})^{(n+1)}$ for example. Otherwise it would not be possible to prove the existence of partial derivatives. In the paper [3], the theorem of Hadamard [4] is used to estimate mixed products of the vectors $\vec{x}_j$. There is a mistake/misprint in the last inequality in [3] which is easy to fix, so the proof is considered to be correct. We are not showing this proof here for it is somewhat cumbersome and bulky one.

Again, we ask if the Cauchy criterion is easy to exploit. In practice, it is absolutely impossible to use the Cauchy criterion to verify differentiability of functions. In this paper, we suggest a simple and easy-to-use criterion for differentiability: that is the Cauchy-like criterion.

## *Cauchy-like criterion for differentiability*

We have seen above, in the case of Geo-criterion, that some necessary condition can be explicitly stated in the formulation of a criterion. Similarly, in the formulation of the Cauchy-like criterion, we require existence of all the partial derivatives. Indeed, if any one of them does not exist then the function is not differentiable anyway. Note that such criterion only works for dimensions no less than two since, in the case of just one variable, the existence of a partial derivative is equivalent to the definition (1).

At the beginning of the paper, we define the partial values of $f(\vec{x})$ as: $f^i = f^i(\vec{x}) = f(x^i \cdot \vec{e}_i) = f(\{0, \ldots x^i, \ldots 0\})$, where $\vec{e}_i$ is the unit vector for the i-th coordinate of $\mathbb{R}^n$. Now, the formulation of the Cauchy-like criterion is as follows.

*Cauchy-like criterion*: For a function $f(\vec{x})$ ($\vec{x} \in \mathbb{R}^n$; $f(0)=0$) to be differentiable at the origin, the following two conditions are necessary and sufficient:
  a) all the partial derivatives exist at the origin; and (5a)
  b) $f(\vec{x}) - \Sigma_i f^i(\vec{x}) = o(\varrho)$, when $\varrho \to 0$, (5b)
     where $\Sigma_i f^i$ – the sum of all partial values of $f(\vec{x})$, and $\varrho = |\vec{x}|$.

This criterion is very easy to use. Indeed, if the condition (5b) does not hold then there is no need to check up on the partial derivatives at all. For example, the function (2) is not differentiable at the origin since the condition (5b) obviously does not hold there. The name "Cauchy-like" has been chosen because partial derivatives are not included in the equation (5b). Let us prove this criterion now.

Proof. *Necessity*:

Let us estimate the left side of the equation (5b) when the condition (1) holds. The vector $\vec{A}$ in the equation (1) is finite, and the components of $\vec{A}$ comprise the set partial derivatives at the origin – the condition (5a). Consequently, for any $1 \leq i \leq n$, $|A^i x^i - f^i| = o(x^i) = o(\varrho)$.

$$|f(\vec{x}) - \Sigma_i f^i| = |f(\vec{x}) - \vec{A}\cdot\vec{x} + \vec{A}\cdot\vec{x} - \Sigma_i f^i| \leq |f(\vec{x}) - \vec{A}\cdot\vec{x}| + |\vec{A}\cdot\vec{x} - \Sigma_i f^i| \leq o(\varrho) + \Sigma_i |A^i x^i - f^i| = o(\varrho).$$

Thus, the necessity of the conditions (5a) and (5b) above is proven.

Proof. *Sufficiency*:

Since all the partial derivatives exist at the origin then their values comprise a vector $\vec{A}$ such that for any $1 \leq i \leq n$, $|A^i x^i - f^i| = o(x^i) = o(\varrho)$. Therefore,

$$|f(\vec{x}) - \vec{A}\cdot\vec{x}| = |f(\vec{x}) - \Sigma_i f^i + \Sigma_i f^i - \vec{A}\cdot\vec{x}| \leq |f(\vec{x}) - \Sigma_i f^i| + \Sigma_i |A^i x^i - f^i| = o(\varrho) + \Sigma_i |A^i x^i - f^i| = o(\varrho).$$

Thus, the sufficiency of the conditions (5a) and (5b) is proven.

Now, let us talk about some other sufficient conditions for differentiability.

As we have seen in *Introduction* (example with the Dirichlet function), there is no criterion for differentiability if one tries to set any conditions on the quality of the function in the vicinity of the origin. The same is true for partial derivatives and their quality. Let us consider the following example of a function $h(\vec{x})$:

$h(\vec{x}) = \varrho^2 \cos(1/\varrho)$ if $\varrho \neq 0$; and $h(\vec{x}=0) = 0$.

According to the Cauchy-like criterion, this function is differentiable at the origin: a) all the partial derivatives exist at the origin, and b) the left side of (5a) for the function $h(\vec{x})$ is equal to $\mathcal{O}(\varrho^2)$. However, any partial derivative of $h(\vec{x})$ is discontinuous at the origin. Consequently, any differentiability conditions that require some quality of the partial derivatives can only be just sufficient conditions. In the next section, we use Cauchy-like criterion to prove some loose sufficient conditions of differentiability.

## *Relaxed sufficient conditions for differentiability*

There are known sufficient conditions for differentiability that require continuity of all the partial derivatives. It turns out that they are too strong. For example, in the case of functions of two variables, it is sufficient to have only one continuous partial derivative. Let us prove these loose or relaxed conditions.

*Relaxed sufficient conditions (2 variables)*: If both partial derivatives of a function f($x,y$) exist at the origin, and one of them is continuous at the origin then f($x,y$) is differentiable at the origin.

Let the "good" variable be $y$. In other words, the partial derivative $f_y(x,y)$ exists in some vicinity of the origin and is continuous at the origin: $f_y(x,y) = f_y(0,0) + \alpha(\varrho)$, where $\alpha(\varrho) \to 0$ when $\varrho \to 0$. The "bad" partial derivative $f_x(0,0)$ simply exists at the origin. According to the Lagrange's mean value theorem, $f(x,y) - f(x,0) = f_y(x,c)y$, where $c = c(x)$, $cy > 0$, while $|c| < |y|$. Since $f_y(x,c)$ is continuous at the origin, we have:

$$f(x,y) - f(x,0) - f(0,y) = f_y(x,c)y - f(0,y) = (f_y(0,0) + \alpha(\varrho))y - f(0,y) = [f_y(0,0)y - f(0,y)] + \alpha(\varrho)y = o(\varrho).$$

Thus, the conditions of the Cauchy-like criterion are met and f($x,y$) is differentiable at the origin. Note here that the "bad" derivative does not participate in the proof at all. It just simply exists to satisfy the condition (5a) of the Cauchy-like criterion. In fact, the "bad" derivative might not even exist in a punctured vicinity of the origin.

Obviously, these relaxed conditions can be generalized for the case of multiple (n) variables:

*Relaxed sufficient conditions (n variables)*: If the first partial derivative of a function f($\vec{x}$) ($\vec{x} \in \mathbb{R}^n$; f(0) = 0) is continuous at the origin over all n variables, and the second partial derivative is continuous over the next (n–1) variables excluding the first one, and so forth, and the next to last derivative is continuous at the origin over the last two variables, while the last partial derivative simply exists at the origin, − then the function f($\vec{x}$) is differentiable at the origin.

These relaxed conditions might seem awkward but they are really loose. In fact, all the textbooks' sufficient conditions for differentiability require continuity of all the partial derivatives over all variables. We place the proof for the general case in the Appendix B at the end of the paper.

So far, we have been considering differentiability of functions in Euclidian space or any normed vector space of finite dimensions. If the domain of the function is not really a vector space then the story is somewhat different.

## *Generalization of the Cauchy-like creterion in non-vector spaces*

In the previous sections, we have discussed differentiability of common real multivariable functions. That is when domain of a function is a normed vector space, also complete – i.e. Banach space, and the codomain is real one-dimensional vector space $\mathbb{R}$. In many cases of applied mathematics, systems analysis and optimization methods, functions are defined in multivariable sets that are not vector spaces. Moreover, the codomain of a function might be itself multi-dimensional linear space. In such cases, the generalization of the differentiability definition (1) is not straightforward and may even be impossible.

Let us consider a finite set of *n* normed vector spaces $\{\mathbb{Y}_i\}$ of finite dimensions. The vector spaces $\mathbb{Y}_i$ may be over different scalar fields and some of them may not even be complete spaces. Let the direct product of $\mathbb{Y}_i$: $\mathbb{P} = \mathbb{Y}_1 \otimes \ldots \mathbb{Y}_i \otimes \ldots \mathbb{Y}_n$, be a domain of a multi-dimensional function F: $\mathbb{P} \to \mathbb{Z}$, where $\mathbb{Z}$ is a codomain of F and a multi-dimensional normed linear space. Thus F is a multivariable function F($\{\mathbf{Y}_i\}$) that maps $\mathbb{P}$ on to $\mathbb{Z}$. Note that $\mathbb{P}$ is a linear space only in the following sense: it includes a zero element $\{0\}$ which plays a role of the origin, and it is linear over additions, that is if any $u, v \in \mathbb{P}$ then $(u \pm v) \in \mathbb{P}$. Thus $\mathbb{P}$ is not actually a vector space because of the different scalar fields of various $\mathbb{Y}_i$. The linear space $\mathbb{P}$ can have the norm assigned in some different ways but they are not essential to our consideration. Also, let $\mathbb{Z}$ be a vector space over a scalar field that includes all the scalar fields of the $\mathbb{Y}_i$ spaces. The norm for the $\mathbb{Z}$ space is defined in a regular way for multi-dimensional vector space over the selected scalar field. All being said now, we can start defining differentiability of the function F($\{\mathbf{Y}_i\}$) analogous to the definition (1).

As before, we consider differentiability at the origin $\{0\}$ and assume that F($\{0\}$) = $\{0\}$ − the zero element in $\mathbb{Z}$ space. First, let us define the partial differentiability of the function F($\{\mathbf{Y}_i\}$). This is straightforward. For any chosen variable $\mathbf{Y}_j$, F must be differentiable at the origin as a function of just only that one variable $\mathbf{Y}_j$. In other words, there must be a linear function $L_j(\mathbf{Y}_j)$: $\mathbb{Y}_j \to \mathbb{Z}$ ($\mathbf{Y}_j \in \mathbb{Y}_j$), so that

$$\|F(\{0, \ldots \mathbf{Y}_j, \ldots 0\}) - L_j(\mathbf{Y}_j)\| = o(\|\mathbf{Y}_j\|). \tag{6}$$

Thus, instead of defining partial derivatives, we define partial differentiability via linear functions $L_j(\mathbf{Y}_j)$ for each variable $\mathbf{Y}_j$.

The definition of differentiability of the function F with domain $\mathbb{P}$ and codomain $\mathbb{Z}$ described above is now as follows. The function F: $\mathbb{P} \to \mathbb{Z}$, is said to be differentiable at the origin {0} if a set of *n* linear functions $L_j(\mathbf{Y}_j)$ exist and the following condition is met in the vicinity of the origin:

$$\|F(\{\mathbf{Y}_1,\ldots \mathbf{Y}_i,\ldots \mathbf{Y}_n\}) - \Sigma_j L_j(\mathbf{Y}_j)\| = o(\varrho), \text{ where } \varrho = \max_j(\|\mathbf{Y}_j\|) \tag{7}$$

The definition (7) is analogous to the definition (1). We see now that $\mathbb{Z}$ needs to be a vector space over a scalar field that entails all the scalar fields of all $\mathbb{Y}_i$ spaces. That is because the codomain of the sum of all linear functions in (7) must coincide with $\mathbb{Z}$. However, since $\mathbb{P}$ is not actually a vector space and may not be even complete space, there is a dramatic difference for differentiability criteria of functions in $\mathbb{P}$ compared to that in $\mathbb{R}^n$. Both Geo-criterion and Cauchy criterion cannot be extended for functions defined in $\mathbb{P}$. Indeed, one cannot define a hyperplane in space that is not any kind of vector space. Also, it is not possible to construct equations analogous to (3) since all the linear functions $L_j(\mathbf{Y}_j)$ are defined in different vector spaces over different scalar fields.

The only allowable extension of differentiability criteria from $\mathbb{R}^n$ onto $\mathbb{P}$ is that of the Cauchy-like criterion. Luckily, the extension is straightforward and simply needs a change of notations as following.

*Cauchy-like criterion*: For a function $F(\{\mathbf{Y}_i\}): \mathbb{P} \to \mathbb{Z}$, (where $\{\mathbf{Y}_i\} \in \mathbb{P}$; $F(\{0\})=\{0\}$) to be differentiable at the origin {0}, the following two conditions are necessary and sufficient:
  a) for any variable $\mathbf{Y}_j$, the function F is partially differentiable; (8a)
  b) $\|F(\{\mathbf{Y}_j\}) - \Sigma_i F^i\| = o(\varrho)$, when $\varrho \to 0$, (8b)
      where $\Sigma_i F^i$ – is the sum of all partial values of $F(\{\mathbf{Y}_i\})$, and $\varrho = \max_j(\|\mathbf{Y}_j\|)$.

The actual proof is exactly the same as before with just minor changes of some notations. Again, the Cauchy-like criterion is very simple and easy to use.

Relaxed sufficient conditions for differentiability can be also formulated in this new settings for $F(\{\mathbf{Y}_i\}): \mathbb{P} \to \mathbb{Z}$. We have proved them for functions in $\mathbb{R}^n$ using Cauchy-like criterion. Apparently, the same can be done for functions in the non-vector space $\mathbb{P}$. For that, we may use the natural definition of continuous partial differentiability.

The partial differentiability of $F(\{\mathbf{Y}_i\}): \mathbb{P} \to \mathbb{Z}$, with respect to the variable $\mathbf{Y}_j$ is said to be continuous at the origin $\{0\} \in \mathbb{P}$, if the equation (6) holds in some small vicinity $U \subset \mathbb{P}$ of the origin, and the function $L_j(\{\mathbf{Y}_i\})$, being a linear function of $\mathbf{Y}_j$ is continuous at the origin. In other words, if $\forall \{\mathbf{Y}_i\} \in U \subset \mathbb{P}$ partial differentiability in U provides $\|F(\mathbf{Y}_1, \ldots \mathbf{Y}_j, \ldots \mathbf{Y}_n) - L_j(\mathbf{Y}_1, \ldots \mathbf{Y}_j, \ldots \mathbf{Y}_n)\| = o(\|\mathbf{Y}_j\|)$ while $\|L_j(\mathbf{Y}_1, \ldots \mathbf{Y}_j, \ldots \mathbf{Y}_n)\| = \mathcal{O}(\|\mathbf{Y}_j\|)$ as a linear function of $\mathbf{Y}_j$, and $\|L_j(\mathbf{Y}_1, \ldots \mathbf{Y}_j, \ldots \mathbf{Y}_n) - L_j(0, \ldots \mathbf{Y}_j, \ldots 0)\| \to 0$ when $\|(\mathbf{Y}_1, \ldots \mathbf{Y}_{j-1}, 0, \mathbf{Y}_{j+1}, \ldots \mathbf{Y}_n)\| \to 0$ as $L_j$ being a continuous function at the origin, then:

$$\|L_j(\mathbf{Y}_1, \ldots \mathbf{Y}_j, \ldots \mathbf{Y}_n) - L_j(0, \ldots \mathbf{Y}_j, \ldots 0)\| = o(\varrho) \tag{9}$$

Note that the relation similar to (9) may have been used in the previous section instead of the Lagrange's mean value theorem. Thus, similar relaxed sufficient conditions for differentiability for $F(\{\mathbf{Y}_i\}): \mathbb{P} \to \mathbb{Z}$, can be proven in the same way using Cauchy-like criterion and equation (9).

## *Discussion*

Let us consider specific *n*-dimensional vector space $\mathbb{C}^n$ and complex functions defined in that space as multivariable functions $F(z_1, \ldots z_i, \ldots z_n): \mathbb{C}^n \to \mathbb{C}$. Complex analysis of such functions is discussed in details in the book by B.V. Shabat [5]. The differentiability criterion is formulated there in the following way.

A function of *n* complex variables $F(z_i): \mathbb{C}^n \to \mathbb{C}$ is differentiable ($\mathbb{C}$-differentiable) at $z \in \mathbb{C}^n$ if and only if:
  a) the Cauchy–Riemann conditions for all $z_i$ variables hold at $z$, and (10a)
  b) F is $\mathbb{R}$-differentiable at $z$ (as a function of 2*n* real variables). (10b)

Note now that for this specific vector space, the condition (10a) is the exact requirement of the partial differentiability similar to those of (5a) and (8a). The second requirement (10b) is more complicated than the simple conditions (5b) and (8b) of the Cauchy-like criterion. Thus both, the Cauchy-like criterion and the Shabat's criterion use the same requirement of partial differentiability, only the former is simpler and easier-to-use.

Partial differentiability over all $\mathbb{Y}_i$ vector spaces is a prerequisite for the Cauchy-like criterion. We did not use the finitude of their dimensions in our consideration. That is, if one can define differentiability of functions in spaces of infinite dimensions then the Cauchy-like criterion may be extended on functions in that type of spaces. Study of differentiability criteria in spaces of infinite dimensions may be very interesting for further investigations.

Summarizing, in our paper we have discussed several differentiability criteria for real functions of multiple variables in *n*-dimensional Euclidean space. We have suggested the simple and easy-to-use Cauchy-like criterion for differentiability and proved it. Relaxed sufficient conditions for differentiability that do not require continuity of all partial derivatives have been also suggested. We have discussed generalization of the Cauchy-like criterion for functions $F(\{Y_i\}): \mathbb{P} \to \mathbb{Z}$, where $\mathbb{P}$ is not any kind of vector space. We suggest that the only valid extension of differentiability criteria from $\mathbb{R}^n$ onto $\mathbb{P}$ is that of the Cauchy-like criterion.

## *References*

1. Alberto Guzman, *Derivatives and Integrals of Multivariable Functions*. *Springer Science & Business Media, 2003.*
2. M. Giaquinta, G. Modica, *Mathematical Analysis: An Introduction to Functions of Several Variables*. *Springer, 2009.*
3. A.M. Macbeath: *A Criterion for Differentiability*. *Edinburgh Math. Notes* (1956), 8-11.
4. G. Hardy, J. Littlewood, G. Polya, *Inequalities*. *Cambridge*, *1934*, p. 34.
5. B.V. Shabat, *Introduction to Complex Analysis: Functions of Several Variables*. *American Mathematical Soc.*, *1992*, Vol. 110.

## *Appendices*

*Appendix A*: Proof of the geometrical criterion for differentiability.

*Geo-Criterion*: If in (*n+1*)-dimensional Euclidean space $\{z, \mathbb{R}^n\}$, there is a hyperplane that is not orthogonal to $\mathbb{R}^n$ and is tangent to the surface $z = f(\vec{x})$ at the origin ($\vec{x} \in \mathbb{R}^n$; $f(0)=0$) then and only then $f(\vec{x})$ is differentiable at the origin.

Proof. *Necessity*:

Let us consider a hyperplane which is given by the equation (1A) in the (n+1)-dimensional Euclidean space $\{z, \mathbb{R}^n\}$ when the statement (1) holds. Let us prove that the hyperplane (1A) is a tangent plane for the surface $z = f(\vec{x})$ at the origin ($\vec{x}=0$).

$$z = \vec{A} \cdot \vec{x} \qquad (1A)$$

If $|\vec{A}|=0$ then the hyperplane (1A) coincides with $\mathbb{R}^n$. Consequently, the distance from the point $\{f(\vec{x}), \vec{x}\}$ to the hyperplane (1A) is equal to $|f(\vec{x})|$. Thus, the condition (1) is essentially the definition of a tangent plane.

If $|\vec{A}| \neq 0$ then the hyperplane (1A) does not coincide with $\mathbb{R}^n$. A vector $\vec{N}_p$ that is normal to the hyperplane may be noted in the following form: $\vec{N}_p = \vec{Z}_e - \vec{A}$, where $\vec{Z}_e$ is the unit vector along the coordinate axis z. Indeed, for any vector $\vec{R}=\{z,\vec{x}\}$ that belongs to the hyperplane (1A), the orthogonality condition holds: $\vec{N}_p \cdot \vec{R} = 0$. The angle between the hyperplanes $\mathbb{R}^n$ and (1A), the angle $\alpha$ between $\vec{N}_p$ and $\vec{Z}_e$, is given by: $\cos\alpha = |\vec{N}_p|^{-1} = (1+|\vec{A}|^2)^{-1/2}$. Since vector $\vec{A}$ has finite components and $\mathbb{R}^n$ does not coincide with (1A), then $0 < \alpha < \pi/2$ and consequently, $\mathbb{R}^n$ is not orthogonal to the hyperplane (1A).

The distance between the point $\{f(\vec{x}), \vec{x}\}$ and the hyperplane (1A) is equal to $d(\vec{x}) = |f(\vec{x}) - \vec{A} \cdot \vec{x}| \cos\alpha$. Thus, from the condition (1), we get that $d(\vec{x}) = o(\varrho)$ with $\varrho = |\vec{x}|$. The distance $r(\vec{x})$ between the origin and the projection of the point $\{f(\vec{x}), \vec{x}\}$ onto the hyperplane (1A) is determined from the right triangle with the hypotenuse $s=|\{f(\vec{x}),\vec{x}\}|$ and the cathetus $d(\vec{x})$: $r = s\{1-(d/s)^2\}^{1/2}$. Note now that if $|\vec{x}| \neq 0$ then $s \neq 0$ and $s \geq \varrho$. Since $d(\vec{x}) = o(\varrho)$ then a small vicinity of the origin in $\mathbb{R}^n$ can be chosen so that $(d/s) \leq 1/2$ for all $\vec{x}$ from that vicinity, and then $r(\vec{x}) \geq (\sqrt{3}/2)s \geq (\sqrt{3}/2)\varrho$. Consequently, $\varrho = \mathcal{O}(r)$ which means that $d(\vec{x}) = o(r)$ and the hyperplane (1A) is tangent to the surface $z = f(\vec{x})$ at the origin.

Proof. *Sufficiency*:

Let a n-dimensional hyperplane be tangent at the origin to the surface $z = f(\vec{x})$ in the (n+1)-dimensional space $\{z, \mathbb{R}^n\}$. This hyperplane is given by some linear functional $L(\{z, \vec{x}\})=0$. The non-orthogonality condition (z axis Å this hyperplane) requires that the coefficient in front of z-coordinate in the $L(\{z, \vec{\mathbf{x}}\})$ is not equal to zero. Hence, the equation $L(\{z, \vec{\mathbf{x}}\})=0$ can be rewritten in a form of equation (1A) with the vector $\vec{A}$ having all components finite.

If this hyperplane coincides with $\mathbb{R}^n$ ($|\vec{A}|=0$) then the tangency condition would be equal to the definition of differentiability (1).

If $|\vec{A}|\neq 0$ then $|f(\vec{\mathbf{x}})-\vec{A}\cdot\vec{\mathbf{x}}|=d(\vec{\mathbf{x}})/\cos\alpha = o(r(\vec{\mathbf{x}}))$ where $r(\vec{\mathbf{x}})$ is the distance between the origin and the projection of the point $\{f(\vec{\mathbf{x}}),\vec{\mathbf{x}}\}$ onto the hyperplane (1A), and $\alpha$ is the angle defined above. When $|\vec{\mathbf{x}}|\neq 0$, $r(\vec{\mathbf{x}})\leq(\varrho^2+|\vec{A}\cdot\vec{\mathbf{x}}|^2)^{1/2}+d\,\mathrm{tg}\,\alpha$. The last inequality is obtained from the triangle in the hyperplane (1A) made by the origin, the projection of the point $\{f(\vec{\mathbf{x}}),\vec{\mathbf{x}}\}$ onto the hyperplane (1A) and the point $\{\vec{A}\cdot\vec{\mathbf{x}},\vec{\mathbf{x}}\}$. Since $d=o(r)$ due to the tangency condition, the immediate consequence of this inequality is that $r=\mathcal{O}(\varrho)$ and $|f(\vec{\mathbf{x}})-\vec{A}\cdot\vec{\mathbf{x}}|=o(r)=o(\varrho)$.

Thus, the Geo-criterion is proven.

*Appendix B*: The general case of loose sufficient conditions for differentiability is considered below.

*Relaxed sufficient conditions (n variables)*: If the first partial derivative of a function $f(\vec{\mathbf{x}})$ ($\vec{\mathbf{x}}\in\mathbb{R}^n$; $f(0)=0$) is continuous at the origin over all n variables, and the second partial derivative is continuous over the next (n-1) variables excluding the first one, and so forth, and the next to last derivative is continuous over the last two variables, while the last partial derivative simply exists at the origin, – then the function $f(\vec{\mathbf{x}})$ is differentiable at the origin.

Let us first prove the following lemma.

*Lemma*: Let the function $f(\vec{\mathbf{x}})$ of n variables ($\vec{\mathbf{x}}\in\mathbb{R}^n$; $f(0)=0$) be such that the function $g(\vec{\mathbf{x}})=f(x^1,\ldots,x^{j-1},0,x^{j+1},\ldots,x^n)$ $=f(\vec{\mathbf{x}}-x^j\vec{\mathbf{e}}_j)$ is differentiable at the origin ($\vec{\mathbf{e}}_j$ is the unit vector of j-th coordinate of $\mathbb{R}^n$). Now, if the j-th partial derivative $f_{x^j}=\partial f/\partial x^j$ is continuous at the origin then $f(\vec{\mathbf{x}})$ is differentiable at $\vec{\mathbf{x}}=0$.

The proof of lemma is analogous to the proof used for two variables.

$f(\vec{\mathbf{x}})-\Sigma_i f^i = f(\vec{\mathbf{x}})-g(\vec{\mathbf{x}})+g(\vec{\mathbf{x}})-\Sigma_i f^i = [f_{x^j}(x^1,\ldots,x^{j-1},c,x^{j+1},\ldots,x^n)x^j-f^j]+[g(\vec{\mathbf{x}})-\Sigma_i f^i+f^j] = [f_{x^j}(0)x^j-f^j]+\alpha(\varrho)x^j+o(\varrho)=o(\varrho)$, here $\alpha(\varrho)\to 0$ when $\varrho\to 0$.

The Cauchy-like criterion holds, and the lemma is proven.

Now, the proof of the general case of the relaxed sufficient conditions for differentiability is by induction. It is true for n=2 variables. Suppose it is true for n−1 variables (n>2). Then, by the lemma proven above it is true for n variables. And the sufficient conditions are proven for the general case.